\documentclass[10pt]{article}
\usepackage{amsthm}
\usepackage{amsmath}
\usepackage{amsfonts}
\usepackage{url}
\usepackage{color,graphics}

\theoremstyle{definition}
\newtheorem{theorem}{Theorem}

\newtheorem{corollary}{Corollary}

\newtheorem{proposition}{Proposition}

\newcommand{\ba}{\begin{array}}
\newcommand{\ea}{\end{array}}
\newcommand{\beq}{\begin{equation}}
\newcommand{\eeq}{\end{equation}}

\def\12{\frac 12}
\def\f2{\frac{r}{2}}

\newcommand{\diag}{\mathrm{diag}\,}
\newcommand{\e}{\epsilon}
\newcommand{\fset}{\mathcal{M}}
\newcommand{\as}{\mathcal{A}}

\begin{document}

\title{A note on the stability number of an orthogonality graph}
\author{E. de Klerk\thanks{Tilburg University. E-mail: e.deklerk@uvt.nl.
Supported by  the Netherlands Organisation for Scientific Research
grant NWO 613.000.214 as well as the NSERC grant 283331 - 04.
Part of this research was performed while on leave from the
Department of Combinatorics and Optimization, University of
Waterloo.} \and D.V. Pasechnik\thanks{Tilburg University. E-mail:
d.v.pasechnik@uvt.nl. {\em Corresponding author}.}} \maketitle

\begin{abstract} 
We consider the orthogonality graph $\Omega(n)$ with $2^n$
vertices corresponding to the vectors $\{0,1\}^n$, two vertices adjacent if
and only if the Hamming distance between them is $n/2$. We show that, for
$n=16$, the stability number of $\Omega(n)$ is $\alpha(\Omega(16))= 2304$, thus
proving a conjecture by Galliard \cite{Galliard thesis}.  The main tool we
employ is a recent semidefinite programming relaxation for minimal distance
binary codes due to Schrijver \cite{Schrijver04}.  

As well, we give a general condition for Delsarte bound on the
(co)cli\-ques in graphs of relations of association schemes to coincide 
with the ratio bound, and use it to show that for $\Omega(n)$ the latter
two bounds are equal to ${2^n}/{n}$.
\end{abstract}


\vspace{0.5cm}

\section{Introduction}
\subsection*{The graph $\Omega(n)$ and its properties}
Let $\Omega(n)$ be the graph on $2^n$ vertices corresponding to
the vectors $\{0,1\}^n$, such that two vertices are adjacent if
and only if the Hamming distance between them is $n/2$. Note
 that $\Omega(n)$ is $k$-regular, where 
$k = \binom{n}{\12 n}$.

It is known that $\Omega(n)$ is bipartite if $n = 2 \mod 4$, and empty if
$n$ is odd.  We will therefore assume throughout that $n$ is a multiple of
$4$. The graph owns its name to another description, in terms of $\pm
1$-vectors.  Then the orthogonality of vectors corresponds to the Hamming
distance $n/2$. 

Moreover, $\Omega(n)$ consists of two isomorphic connected 
components, $\Omega_0(n)$, respectively $\Omega_1(n)$,
containing all the vertices of even, respectively odd, 
Hamming weight. 
For a detailed discussion of the
properties of $\Omega(n)$, see Godsil \cite{Godsil2004}, the PhD thesis
of Newman \cite{thesisnewman}, and \cite{GodsilNewman2005}.

\begin{figure}[!ht]
\begin{center}
\input{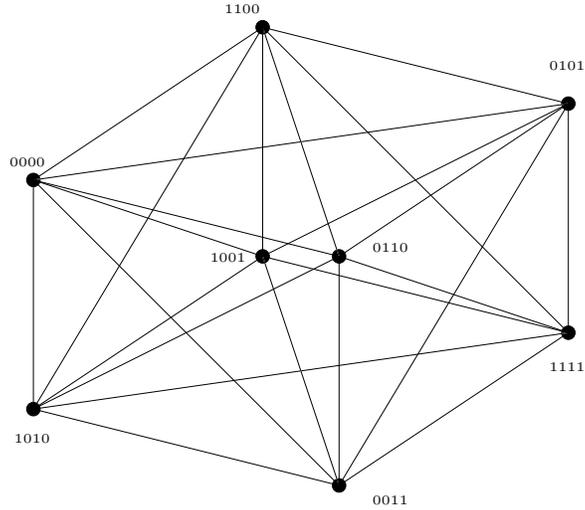}
\caption{The connected component of $\Omega(4)$ corresponding to vertices of
even Hamming weight.}
\end{center}
\end{figure}

In this note we study upper bounds on the stability number $\alpha(\Omega(n))$.

Galliard \cite{Galliard thesis} pointed out the following way of
constructing maximal stable sets in $\Omega(n)$.  Consider the component
$\Omega_{\e}(n)$ of $\Omega(n)$, 
for $1-\e=\frac{n}{4}\mod 2$, and take all
vertices of Hamming weight $\e,\e+2,\ldots,\e+2\ell,\ldots,n/4-1$. 
Obviously, these
vertices form a stable set $S$ of $\Omega(n)$ of size
\[
\sum_{i=\e}^{ \lfloor n/8\rfloor} \binom{n}{2i-\e}.
\]
We can double the size of $S$ by adding the bit-wise complements of
the vertices
in $S$, and double it again by taking the union with the corresponding
stable set in $\Omega_{1-\e}(n)$.
Thus we find that
\beq
\label{lower bound}
\alpha(\Omega(n)) \ge 4
\sum_{i=\e}^{ \lfloor n/8\rfloor} \binom{n}{2i-\e}
 := {\underline \alpha}(n)
\eeq
For $n= 16$ this evaluates to $\alpha(\Omega(n)) \ge 2304$.
Galliard et al \cite{Galliard2002} could show that
$\alpha(\Omega(16)) \le 3912$.
In this note we will show that, in fact,  $\alpha(\Omega(16)) = 2304$.
This was conjectured by Galliard \cite{Galliard thesis}, and
Newman \cite{thesisnewman} has recently conjectured
that the value (\ref{lower bound}) actually equals $\alpha(\Omega(n))$ whenever $n$ is a multiple of $4$.

\subsection*{A quantum information game}
One motivation for studying the graph $\Omega(n)$
comes from quantum information theory.
Consider the following game from \cite{Galliard2002}.

Let $r\ge 1$ and $n = 2^r$. Two players, A and B, are asked the questions $x_A$ and $x_B$,
coded as $n$-bit strings satisfying
\[
d_H(x_A,x_B) \in \left\{0,\12n\right\}
\]
where $d_H$ denotes the Hamming distance.
A and B win the game if they  give answers $y_A$ and $y_B$,
coded as binary strings of length $r$ such that
\[
y_A = y_B \Longleftrightarrow x_A = x_B.
\]
A and B are not allowed any communication (except a priori deliberation).

It is known that A and B can always win the game if their $r$ output bits are {\em maximally
entangled quantum bits} \cite{Brassard99} (see also \cite{thesisnewman}).

For classical bits, it was shown by Galliard et al \cite{Galliard2002} that the game cannot always be won
if $r = 4$.
The authors proved this by pointing out that whether or not the game can always be won
is equivalent to the question
\[
\chi(\Omega(n)) \le n ?
\]
Indeed, if $\chi(\Omega(n)) \le n$ then A and B may color $\Omega(n)$  a priori using $n$ colors.
The questions $x_A$ and $x_B$ may then be viewed as two vertices of $\Omega(n)$,
and the A and B may answer their respective questions by giving the color of the vertices $x_A$ and $x_B$ respectively,
coded as binary strings of length $\log_2 n = r$.

Galliard et al.\ \cite{Galliard2002} showed that $\chi(\Omega(16)) > 16$, i.e. that the game cannot be won for $n=16$.
They proved this by showing that $\alpha(\Omega(16)) \le 3912$ which
implies
\[
\chi(\Omega(16)) \ge \left\lceil \frac{2^{16}}{\alpha(\Omega(16))} \right\rceil
 \ge \left\lceil \frac{2^{16}}{3912} \right\rceil  = 17.
\]
In this note we sharpen their bound by showing that $\alpha(\Omega(16)) = 2304$, which implies $\chi(\Omega(16)) \ge 29$.

Our main tool will be a semidefinite programming bound on $\alpha(\Omega(n))$
that is due to Schrijver \cite{Schrijver04}, where it is formulated
for minimal distance binary codes.

\section{Upper bounds on $\alpha(\Omega(n))$}
In this section we give a review of known upper bounds on $\alpha(\Omega(n))$
and their relationship.

\subsection{The ratio bound}
The following discussion is condensed from Godsil \cite{Godsil2004}.

\begin{theorem}
Let $G = (V,E) $ be a $k$-regular graph with adjacency matrix
$A(G)$, and let $\lambda_{\min}(A(G))$ denote the smallest eigenvalue of $A(G)$.
Then
\beq\label{eq:rb}
\alpha(G) \le \frac{|V|}{1 - \frac{k}{\lambda_{\min}(A(G))}}.
\eeq
\end{theorem}
This bound is called the {\em ratio bound}, and was first derived by Delsarte
\cite{delsarte73} for graphs in association schemes (see Sect.~\ref{sect:db}
for more on the latter).

Recall that $\Omega(n)$ is $k$-regular with $k = \binom{n}{\12n}$.
Ignoring multiplicities, the spectrum of $\Omega(n)$ is given by
\beq\label{eq:eigomega}
\lambda_m =
\frac{2^{\12n}}{(\12 n)!}
(m - 1)(m - 3) \cdots (m - n + 1)  \quad\quad (m = 1,\ldots,n).
\eeq
The minimum is reached at $m=2$, and we get
\beq\label{eq:eigminomega}
\lambda_{\min}(A(\Omega(n)) = \frac{2^{\12n}}{(\12 n)!}
(1)(-1)(-3) \cdots ( -n + 3) = -\frac{\binom{n}{\12n}}{n-1}.
\eeq
The ratio bound therefore becomes
\beq
\label{ratio bound}
\alpha(\Omega(n)) \le \frac{2^n}{n}.
\eeq
This is the best known upper bound on $\alpha(\Omega(n))$, but it is known that this bound is not tight:
Frankl and R\"odl  \cite{Frankl} showed that there exists some $\e > 0$ such that $\alpha(\Omega(n)) \le {(2-\e)^n}$.
For specific (small) values of $n$ one can improve on the bound (\ref{ratio bound}), as we will show for $n \le 32$.

\subsection{The Delsarte bound and $\vartheta'$}\label{sect:db}
Here we are going to use more linear algebra that naturally arise around
$\Omega(n)$. We recall the following definitions,
cf. e.g. Bannai and Ito \cite{BanIto84}.

\paragraph{Association schemes.} 
An association scheme $\as$ is a commutative subalgebra
of the full $v\times v$-matrix algebra with a distinguished basis
$(A_0=I,A_1,\dots,A_n)$ of 0-1 matrices, with an extra property that
$\sum_i A_i$ equals the all-ones matrix. 
One often views $A_j$, $j\geq 1$,  as the adjacency matrix of 
a graph on $v$ vertices;  $A_j$ is often referred to as the $j$-th
{\em relation} of $\as$.
As the $A_j$'s commute,
they have $n+1$ common eigenspaces $V_i$. Then $\as$ is isomorphic,
as an algebra, to the algebra of diagonal matrices 
$\diag(P_{0,j},\dots,P_{n,j})$, where $P_{ij}$ denotes the
eigenvalue of $A_j$ on $V_i$. The matrix $P=(P_{ij})$ is called
{\em first eigenvalue matrix} of $\as$. 
The set of $A_j$'s is closed under taking
transpositions: for each $0\leq j\leq n$ 
there exists $j'$ so that $A_j=A_{j'}^T$. In particular, 
$P_{ij}=\overline{P_{ij'}}$. An association scheme with all $A_j$ symmetric
is called {\em symmetric}, and here we shall consider such schemes only.
There is a 
matrix $Q$ (called {\em second eigenvalue matrix}) 
satisfying $PQ=QP=vI$.
In what follows it is assumed (as is customary in the
literature) that the eigenspace $V_0$ corresponds
to the eigenvector $(1,\dots,1)$; then the 0-th row of $P$
consists of the degrees $v_j$ of the graphs $A_j$.
It is remarkable that the 0-th row of $Q$ consists of dimensions
of $V_i$. 
\medskip

Let
$\vartheta'$ denote the Schrijver $\vartheta'$-function \cite{Schr79}:
\[
\vartheta'(G) = \max \left\{\mbox{Tr} \left(JX\right) \; : \; \mbox{Tr}\left(AX\right)  =
0,\; \mbox{Tr} (X)  =  1, \;  X  \succeq 0, \; X \ge 0
\right\}.
\]
For any graph $G$ one has $\alpha(G) \le \vartheta'(G)$. Moreover,
$\vartheta'(G)$
is smaller than or equal to the ratio bound (\ref{eq:rb}) for regular graphs, 
as noted by Godsil \cite[Sect.~3.7]{Godsil2004}.

For graphs with adjacency matrices of the form 
$\sum_{j\in\fset} A_j$, with $\fset\subset\{1,\dots,n\}$ and
$A_j$'s from  the 0-1 basis of an
association scheme $\as$,
the bound $\vartheta'$ coincides, as was proved by Schrijver \cite{Schr79},
 with the following bound due to Delsarte
\cite{delsarte72,delsarte73} 
\begin{equation}\label{eq:delslp}
\max 1^T w\ \text{subject to }
w\geq 0,\ Q^T w\geq 0,\ w_0 = 1, w_j = 0\ \text{for }j\in\fset,
\end{equation}
where $Q$ is the second eigenvalue matrix of $\as$.

The bound (\ref{eq:delslp}) is often stated for 
(and was originally developed for) bounding the maximal size
of a $q$-ary code of length $n$ and minimal distance $d$; 
then the association scheme $\as$ becomes the Hamming distance
association scheme $H(n,q)$ and $\fset=\{1,\dots,d-1\}$.
The relations of $H(n,q)$ can be viewed as graphs on the 
vertex set of $n$-strings on $\{0,\dots,q-1\}$:  
the $j$-th graph of $H(n,q)$ is given by
\[
(A_j)_{XY}=\left\{ \ba{ll}  
1 & \text{if  }d_H(X,Y)=j\\
0 & \text{otherwise}
\ea \right.
\]

For $H(n,q)$ the first and the second eigenvalue
matrices $P$ and $Q$ coincide, and are given by $P_{ij}=K_j(i)$,
where $K_k$ is the {\em Krawtchouk polynomial}
\[
K_k(x) := \sum_{j=0}^k (-1)^j(q-1)^{k-j}\binom{x}{j}\binom{n-x}{k-j}.
\]

For $\Omega(n)$, the bound (\ref{eq:delslp}) is as above with $\as=H(n,2)$
and $\fset=\{\frac{n}{2}\}$.  Newman \cite{thesisnewman} has shown
computationally that $ \vartheta'(\Omega(n)) = 2^n/n \mbox{ if $n \le
64$}, $ i.e.\ the ratio and $\vartheta'$ bounds coincide for $\Omega(n)$
if $n \le 64$.  We show that it is the case for all $n$, as an easy
consequence of the following.  
\begin{proposition}\label{prop:delseqrat}
Let $\as$ be an association scheme with the 0-1 basis $(A_0,\dots,A_n)$
and eigenvalue matrices $P$ and $Q$.  Let $A_r$ have the least eigenvalue
$\tau=P_{\ell r}$ and assume \[ v_r P_{\ell i}\geq v_i\tau,\quad 0\leq
i\leq n.  \] Then the Delsarte bound (\ref{eq:delslp}), with
$\fset=\{r\}$, and the ratio bound (\ref{eq:rb}) for $A_r$ coincide.
\end{proposition} \begin{proof} Let $P_j$ denote the $j$-th row of $P$. 

As we already mentioned, the bound (\ref{eq:rb}) for regular graphs
always majorates  (\ref{eq:delslp}). 
Thus it suffices to present a feasible vector for 
the LP in (\ref{eq:delslp}) that gives the objective value the
same as (\ref{eq:rb}). 

We claim that
$$a=\frac{-\tau}{v_r-\tau}P_0^T+
\frac{v_r}{v_r-\tau}P_\ell^T$$ 
is such a vector. It is straightforward to check that 
$a_0=1$ and $a_r=0$, as required. 
By the assumption of the proposition, $a\geq 0$.
As $PQ=vI$, any nonnegative 
linear combination $z$ of
the rows of $P$ satisfies $Q^T z^T\geq 0$. As 
$a^T$ is such a combination, we obtain $Q^T a\geq 0$.

Finally, to compute $1^T a$, note that $1^T P_0^T=v$ and 
$1^T P_\ell^T=0$.
\end{proof}

\begin{corollary}
The bounds (\ref{eq:delslp}) and (\ref{eq:rb}) coincide for $\Omega(n)$. 
\end{corollary}
\begin{proof}
We apply Proposition~\ref{prop:delseqrat} to $\as=H(n,2)$ and 
$r=\frac{n}{2}$. 
Then the eigenvalues of  $A_{r}=\Omega(n)$ given
in (\ref{eq:eigomega})
comprise the $r$-th column on $P$, in particular
the least eigenvalue $\tau$ 
equals $P_{2,r}$, by (\ref{eq:eigminomega}) above. 
The assumption of the proposition translates 
into\footnote{Here $m!!=m(m-2)(m-4)\dots$, the {\em double factorial}.}
\[
\binom{n}{\frac{n}{2}}K_i(2)-\binom{n}{i}K_{\frac{n}{2}}(2)=
\frac{2^{\frac{n}{2}+2} (n-2)! (n-1)!! (\frac{n}{2}-i)^2}
 {i! (\frac{n}{2})! (n-i)!}\geq 0,
\]
as claimed.
\end{proof}

\subsection{Schrijver's improved SDP-based bound}
Recently, Schrijver \cite{Schrijver04} has suggested a new SDP-based bound for minimal distance codes, that
is at least as good as the $\vartheta'$ bound, and still of size polynomial in $n$.
It is given as the optimal value of a semidefinite programming (SDP) problem.

In order to introduce this bound (as applied to $\alpha(\Omega(n))$) we require some notation.

For $i,j,t \in \{0,1,\ldots,n\}$, and  $X,Y \in \{0,1\}^n$  
define the matrices
\[
\left( M^t_{i,j}\right)_{X,Y} = \left\{
\ba{ll}
1 & \mbox{if $|X| = i$, $|Y| = j$, $d_H(X,Y) = i+j-2t$} \\
0 & \mbox{otherwise}
\ea \right.
\]
The upper bound is given as the optimal value of the following semidefinite program:
\[
\bar \alpha(n) := \max \sum_{i=0}^n \binom{n}{i} x_{i,0}^0
\]
subject to
\begin{eqnarray*}
x_{0,0}^0 &=& 1 \\
0 &\le & x_{i,j}^t \le x_{i,0}^0 \mbox{ for all $i,j,t \in \{0,\ldots,n\}$} \\
x_{i,j}^t &= &x_{i',j'}^{t'} \mbox{ if $\{i',j',i'+j'-2t'\}$ is a permutation of  $\{i,j,i+j-2t\}$} \\
x_{i,j}^t &= & 0 \mbox{ if $\{i,j,i+j-2t\} \cap \{\12n\} \neq \emptyset,$}
\end{eqnarray*}
as well as
\[
\sum_{i,j,t} x_{i,j}^tM_{i,j}^t \succeq 0, \;\; \sum_{i,j,t} \left(x_{i+j-2t,0}^0-x_{i,j}^t\right)M_{i,j}^t \succeq 0.
\]
The matrices $ M^t_{i,j}$ are of order $2^n$ and therefore too large to
compute with in general.
Schrijver pointed out that these matrices form a basis of the Terwilliger
algebra of the Hamming scheme,
and worked out the details for computing the irreducible block diagonalization
of this (non-commutative) matrix algebra of dimension $O(n^3)$.

Thus, analogously to the $\vartheta'$-case, the constraint $\sum_{i,j,t} x_{i,j}^tM_{i,j}^t \succeq 0$ is replaced by
\[
\sum_{i,j,t} x_{i,j}^tQ^TM_{i,j}^tQ \succeq 0
\]
where $Q$ is an orthogonal matrix that gives the irreducible block
diagonalization.
For details the reader is referred to Schrijver \cite{Schrijver04}.
Since SDP solvers can exploit block diagonal structure, this reduces the sizes
of the matrices
in question to the extent that computation is possible in the range $n \le 32$.


\subsection{Laurent's improvement}
In Laurent \cite{Laurent05} one finds a study placing the relaxation
\cite{Schrijver04} into the framework of {\em moment sequences} of
\cite{Lass, Laurent03}. This study also explains the relationship with known lift-and-project methods
for obtaining hierarchies of upper bounds on $\alpha(G)$.

Moreover, Laurent \cite{Laurent05} suggests
a refinement of the Schrijver relaxation that takes the following form:
\[
l_+(n) := \max 2^n  x_{0,0}^0
\]
subject to
\begin{eqnarray*}
0 &\le & x_{i,j}^t \le x_{i,0}^0 \mbox{ for all $i,j,t \in \{0,\ldots,n\}$} \\
x_{i,j}^t &= &x_{i',j'}^{t'} \mbox{ if $\{i',j',i'+j'-2t'\}$ is a permutation of  $\{i,j,i+j-2t\}$} \\
x_{i,j}^t &= & 0 \mbox{ if $\{i,j,i+j-2t\} \cap \{\12n\} \neq \emptyset,$}
\end{eqnarray*}
as well as
\[
\sum_{i,j,t} x_{i,j}^tM_{i,j}^t \succeq 0
\]
and
\[
\left(
\ba{cc}
1- x_{0,0}^0 & c^T \\
c & \sum_{i,j,t} \left(x_{i+j-2t,0}^0-x_{i,j}^t\right)M_{i,j}^t
\ea
\right) \succeq 0,
\]
where $c := \sum_{i=0}^n \left(x_{0,0}^0 - x_{0,i}^0\right)\chi_i$, and
$\chi_i$ is defined by
\[
(\chi_i)_X := \left\{ \ba{ll}
1 & \mbox{if $|X| = i$} \\
0 & \mbox{else.}
\ea
\right.
\]
This SDP problem may be block-diagonalised as before to obtain an SDP of size
$O(n^3)$.

\section{Computational results}
To summarize, the bounds we have mentioned satisfy:
\[
{\underline \alpha}(n) \le \alpha(\Omega(n)) \le l^+(n) \le \bar \alpha(n) \le \vartheta'(\Omega(n)) =2^n/n.
\]

In Table \ref{Table:1} we show the numerical values for $\bar \alpha(n)$
and $l_+(n)$ that were obtained using the SDP solver SeDuMi by Sturm
\cite{SeDuMi}, with Matlab 7 on a Pentium IV machine with 1GB of memory.
Maltab routines that we have written to generate the corresponding
SeDuMi input are available online \cite{input-bounds}.

\begin{table}[!ht] \label{tab:table1}
\begin{center}
\begin{tabular}{|c|c|c|c|c|} \hline
$n$ & ${\underline \alpha}(n)$ &   $l_+(n)$ & $\bar \alpha(n)$  & 
$\vartheta'(\Omega(n))=\lfloor 2^n/n \rfloor$ \\ \hline
$16$& $2304$ & $2304$ & $2304$ & $4096$ \\
$20$& $20,144$ & $20,166.62$  & $20,166.98$ & $52,428$  \\
$24$& $178,208$ & $183,373$ & $184,194$ & $699,050$  \\
$28$ & $1,590,376$ & $1,848,580$ & $1,883,009$ & $9,586,980$ \\
$32$& $14,288,896$ & $21,103,609$&  $21,723,404$ & $134,217,728$ \\ \hline
\end{tabular}
\end{center}
\caption{\label{Table:1} Lower and upper bounds on $\alpha(\Omega(n))$.}
\end{table}
Note that the lower and upper bounds coincide for $n=16$, proving
that $\alpha(\Omega(16)) = 2304$. The best previously known upper
bound, obtained by an {\em ad hoc} method,
was $\alpha(\Omega(16)) \le 3912$ \cite{Galliard2002}.

The value $\bar \alpha(20)= 20,166.98$ implies
that $$\alpha(\Omega(20)) \in \{20144,20148,20152,20156,20160,20164\}$$ since $\alpha(\Omega(n))$
is always a multiple of $4$.
Another implication is that $n=20$ is the smallest value
of $n$ where the upper bounds
$\bar \alpha(n)$ and $l_+(n)$ are not tight.

It is worth noticing that the Schrijver and Laurent bounds ($\bar \alpha(n)$
and $l_+(n)$ respectively) give relatively big improvements over the Delsarte
bound $\frac{2^n}{n}$. This is in contrast to the relatively small improvements
that these bounds give for binary codes, cf.\ \cite{Schrijver04,Laurent05}.  We
also note that these relaxations are numerically ill-conditioned for $n \ge
24$. This makes it difficult to solve the corresponding SDP problems to high
accuracy. The recent study by De Klerk, Pasechnik, and Schrijver \cite{deKPS05}
suggests a different way to solve such SDP problems, leading to larger SDP
instances, but which may avoid the numerical ill-conditioning caused by
performing the irreducible block factorization.

\section*{Acknowledgements} The authors would like to thank Chris Godsil for
communicating this problem to them and for fruitful comments, and 
Willem Haemers, Mike Newman
and Lex Schrijver for useful discussions.

\end{document}